\documentclass[12pt,a4paper]{article}
\usepackage{cmap}
\usepackage[T1]{fontenc}
\usepackage{amsfonts,amsmath,amssymb,amsthm,longtable,hhline}
\usepackage{bm}
\usepackage{cite}
\usepackage{color,soul}
\usepackage{color,xcolor}
\usepackage{multicol}
\usepackage{tabularx,caption,multirow}
\usepackage{graphicx}

\usepackage[labelformat=empty,labelsep = none,skip=0pt,justification=raggedright,singlelinecheck=false]{caption}

\usepackage[linktoc=none,linktocpage=true,unicode,bookmarks,bookmarksopen,bookmarksopenlevel=2,colorlinks,linkcolor=blue,citecolor=blue,urlcolor=blue]{hyperref}

\usepackage{tcolorbox}
\usepackage{multirow}
\usepackage{tikz}
\usepackage{xparse}

\renewcommand{\arraystretch}{1.2} 

\makeatletter
\renewcommand{\@biblabel}[1]{#1.} 
\makeatother
\theoremstyle{remark}

\newcounter{urav}[section]
\newcounter{resh}[urav]

\newcounter{exmp}

\let\ds=\displaystyle
\let\ts=\textstyle

\def\fracskip{\mskip 1mu \relax}
\def\nfrac#1#2{{\fracskip#1\fracskip\over\fracskip#2\fracskip}}
\def\dfrac#1#2{{\ds\nfrac{#1}{#2}}}
\def\tfrac#1#2{{\ts\nfrac{#1}{#2}}}
\let\frac=\nfrac
\def\pd#1#2{\dfrac{\partial#1}{\partial#2}}
\def\pdd#1#2#3{\ifx#2#3\pd{^2#1}{#2^2}\else\pd{^2#1}{#2\partial#3}\fi }

\newcommand{\clh}[1]{\colorbox{yellow}{#1}}%
\newcommand{\clhp}[1]{\colorbox{yellow}{\parbox{\textwidth}{#1}}}%

\newlength\aptextwidth
\definecolor{BrickRed}{rgb}{0.588,0.098,0.055}

\def\noblue#1{\ifmmode \text{#1}\else #1\fi}
\def\noclh#1{\ifmmode \text{#1}\else #1\fi}
\def\rem#1{}

\let\le=\leqslant

\def\ttable#1. #2{\begin{table}[t]\tablehat{#1}{#2}}
\def\mtable#1. #2{\begin{table}[hbtp]\tablehat{#1}{#2}}
\def\ptable#1. #2{\begin{table}[p]\tablehat{#1}{#2}}
\def\tablehat#1#2{\centering\small \vbox{\parindent=0pt
  \leftskip=0pt plus.5\hsize \rightskip=\leftskip \parfillskip=0pt
  ТАБЛИЦА #1\\ #2}\nobreak\medskip\medskip }
\def\texendtable{\end{table}}
\def\notation{\par\ifnum\lastpenalty<25000 \bigbreak \fi
  \noindent\triangle\enspace\ignorespaces}
\def\ccbox#1{$\vcenter{\ialign{\hfil##\hfil\cr#1\crcr}}$}

\def\lcbox#1{$\vcenter{\ialign{##\hfil\cr#1\crcr}}$}

\def\hline{\noalign{\hrule}}
\def\rlineskip{\vruleskip\hline\vruleskip}

\let\ds=\displaystyle
\let\ts=\textstyle

\let\rightleftharpoons=\rightleftarrows

\def\eqnitemskip{\ifhmode \else \par
  \ifnum\lastpenalty>24999
    \ifnum\lastpenalty=25004 \fi
  \else \medbreak \fi \fi }
\def\eqnitem #1. {\eqnitemskip
  {\setbox0=\hbox{$#1^\circ$.\enspace}%
  \ifdim\wd0>\parindent \box0\ignorespaces \else
  \hbox to\parindent{\unhbox0\hss}\ignorespaces\fi}}
\def\eqnitemnobreak #1. {\noindent
  {\setbox0=\hbox{$#1^\circ$.\enspace}%
  \ifdim\wd0>\parindent \box0\ignorespaces \else
  \hbox to\parindent{\unhbox0\hss}\ignorespaces\fi}}
\newdimen\eqnparindent
\eqnparindent=\parindent
\def\eqnitem #1. {\eqnitemskip\noindent\hskip\eqnparindent $#1^\circ$.\enspace\ignorespaces }
\def\eqnitemnobreak #1. {\noindent\hskip\eqnparindent $#1^\circ$.\enspace\ignorespaces }
\def\simpleitem #1. {\eqnitemskip\noindent\hskip\eqnparindent #1.\enspace\ignorespaces }

\def\eqalignno#1{\displ@y \tabskip\centering
  \halign to\displaywidth{\hfil$\@lign\displaystyle{##}$\tabskip\z@skip
    &$\@lign\displaystyle{{}##}$\hfil\tabskip\centering
    &\llap{$\@lign\eqnofont##$}\tabskip\z@skip\crcr
    #1\crcr}}
\let\eqalignno=\eqalignm
\def\eqcenter#1{\displ@y \tabskip\centering
  \halign{\hfil$\displaystyle{##}$\hfil\crcr
    #1\crcr}}
\def\eqcenterno#1{\displ@y \tabskip\centering
  \halign to\displaywidth{\hfil$\@lign\displaystyle{##}$\hfil
    \tabskip\centering&\llap{$\@lign\eqnofont##$}\tabskip\z@skip\crcr
    #1\crcr}}
\def\texcases#1{\left\{\,\vcenter{\normalbaselines\m@th
    \ialign{$##\hfil$&\quad##\hfil\crcr#1\crcr}}\right.}
\def\Displaylines#1{\vcenter{\displ@y \tabskip\z@skip
  \halign{\hbox to\displaywidth{$\@lign\hfil\displaystyle##\hfil$}\crcr
    #1\crcr}}}

\usepackage{caption}
\begin{document}

\centerline{\bf\Large  Nonclassical symmetries of polynomial equations} 
\centerline{\bf\Large  and test problems with parameters}
\centerline{\bf\Large  for computer algebra systems\clh{$^*$}}

\bigskip

\centerline{Inna K. Shingareva$^{a}$, Andrei D. Polyanin$^{b}$}
\medskip
\centerline{\it $^a$ Department of Mathematics, University of Sonora,} 
\centerline{\it Blvd. Luis Encinas y Rosales S/N, Hermosillo C.P. 83000, Sonora, M\'exico}
\centerline{\it $^b$ Ishlinsky Institute for Problems in Mechanics RAS,}  
\centerline{\it 101 Vernadsky Avenue, bldg 1, 119526 Moscow, Russia}
\centerline{inna.shingareva@unison.mx, polyanin@ipmnet.ru}
\bigskip
\bigskip


\let\thefootnote\relax\footnotetext{
\hskip-20pt\clhp{$^*$ 
This work is an expanded English version of our article, published in Spanish:
\url{https://sahuarus.unison.mx/index.php/sahuarus/article/view/180}.}}

\begin{abstract}
Nonclassical symmetries and reductions of polynomial  equations and
systems of polynomial equations are considered.
It is shown that specific polynomial equations having ``hidden" symmetries
can be reduced to classical symmetric systems of polynomial equations
by introducing a new additional variable.
It has been established that symmetric systems of polynomial equations of mixed type,
consisting of symmetric and anti-symmetric polynomials, can be transformed into simpler systems.
A method is presented for solving nonclassical symmetric systems
of two polynomial equations that change places when the unknowns are permuted.
We study polynomial equations containing the second iteration of a given polynomial,
which are reduced to nonclassical symmetric systems of equations.
New higher-degree polynomial equations containing free parameters that admit solutions in radicals are found.
Three such equations of the sixth and ninth degrees are further used as test problems with parameters for
analyzing the capabilities of two leading computer algebra systems.
It is shown that currently, the Maple and Mathematica systems do not allow us
to efficiently find analytical solutions (in radicals) of polynomial equations
with free parameters, but they allow us to obtain numerical solutions
of equations for fixed numerical values of the parameters.
The results of this work and the proposed test problems with parameters
can be used  to further improve existing computer algebra systems.
\end{abstract}


\section{Introduction}\label{sec1}

A polynomial (algebraic) equation of $n$\,th degree containing one unknown quantity~$x$
is an equation  of the form
\begin{equation}
P_n(x)=0,
\label{Eq1}
\end{equation}
where $P_n(x)$  is the $n$\,th degree polynomial.

The solution (root) of equation \eqref{Eq1} is the number $x_*$ that satisfies the condition
$P_n(x_*){=}0$.
Solving equation~\eqref{Eq1} means finding all its solutions.

By symmetries of a mathematical equation, we mean transformations
that preserve the form of the equation under consideration.
Reduction is a way of transforming a mathematical equation
into a simpler or more convenient (from some point of view)
form for analysis and solution.

For centuries, the main problem of algebra was
the development of methods for solving polynomial equations.
Formulas for solving polynomial equations of the first and second degrees
have been known since ancient times.

In the 16th century, Italian mathematicians Scipione del Ferro,
Niccol\'o Fontana (a.k.a. Tartaglia), Gerolamo Cardano,
and Ludovico Ferrari obtained formulas for solving polynomial equations
of the third and fourth degree.
Historical information about this can be found in
\cite{Turnbull1947,Waerden1985}
(see also \cite{bro2015,kor2000,polman2007,yac2012,the2016,Sanchez2020,Prodanov2021,cha2023},
where various representations of solutions to these
equations are given).
\medskip

\textit{Remark 1}.
Much earlier, in the 11\,th century, Omar Khayyam, an outstanding poet,
philosopher and scientist, wrote ``Treatise on Demonstration of Problems of Algebra",
where he described a geometric method
for solving polynomial equations of the third degree,
having positive roots, using conic sections (see \cite{Waerden1985,Siadat2021}).
Later, in the 17th century, a theory for solving cubic equations
based on conic sections was developed by Ren\'e Descartes
and by other scientists who were not familiar with the works of Khayyam.
\medskip

Between the 16th and early 19th centuries,
for polynomial equations of degree higher than the fourth,
various attempts have been made to represent
{\it solutions in radicals}, i.e. in the form of an expression
containing only the coefficients of a polynomial equation and
the operations of addition, subtraction, multiplication, division and root extraction.
However, these attempts were unsuccessful in the general case.

In the 19th century, Paolo Ruffini, Niels Hendrik Abel, and \'Evariste Galois established
that in the general case the solution of polynomial equations of degree
higher than the fourth cannot be expressed in radicals
\cite{Waerden1985,Struik1986} (see also \cite{kin1996,bro2015,polman2007,pol2024,Ramond2020}).
Various theorems on the distribution of real and complex roots
of polynomial equations can be found, for example, in \cite{bro2015,polman2007,bha2023}.

It is important to note that the solution of specific polynomial equations
of the fifth degree and higher can be expressed in radicals.
Some such equations can be found in the handbook~\cite{pol2024},
where in addition to solutions in radicals, solutions are also given
that are expressed in terms of special functions
(for polynomial equations of the fifth degree, see also \cite{kin1996,spe1994}).
In several cases, it is possible to find the roots of polynomial equations and systems of polynomial equations
that have different symmetry properties and therefore allow for significant simplifications 
(see, for example, \cite{pol2024,col1997,bol2002,cor2009}).

Nowadays, computer algebra systems (CAS) are being used
for solving polynomial equations and systems of polynomial equations,
and they provide good results (as well as MATLAB and other computational software)
for finding roots of polynomial equations with numerical coefficients,
i.e. without parameters
(see \cite{rou1999,lazard1999,xia2002,rouzim2003,shili2009,petkovic2017,petkovic2020,herzberger2003,petkovic2019,petkovic2022}).
However, as will be shown later in this article,
computer algebra systems are presently ineffective
for finding all solutions to polynomial equations
that contain free parameters.
Nevertheless, the computer algebra system Maple can find some (but not all)
real roots of polynomial equations with parameters
(see, for example, \cite{chen2013,chistov2018,kapur2013,kapur2021,lazard2007,lazard2009,liang2009,montesa2010}
and Section~9 of this work).

It should be noted that it is important to create test problems
for development, verification, and improvement
of computer algebra systems.
For this purpose, based on the analytical results obtained
in the first part of the article, we propose several
test problems with parameters for higher-degree polynomial equations.
We give some examples of nontrivial polynomial equations of the sixth and ninth degrees
(containing free parameters) that admit solutions in radicals.

Finally, in Sections 8 and 9, these test problems for high-degree polynomial equations with parameters
were used to test the most common computer algebra systems Maple and Mathematica~\cite{shili2009}.
\medskip

\textit{Remark 2}.
It should be noted that Maple and Mathematica can be used quite successfully
to find exact solutions to ordinary differential equations (ODEs) with parameters
(see, for example, \cite{shili2018,precup2018,marasco2001}).
\medskip

To formulate test problems with parameters, we first consider
some polynomial equations and systems
of polynomial equations that have certain symmetry properties
and allow various simplifications and solutions in radicals.

\section{Symmetric systems of polynomial equations}\label{ss:1.2.4}

\eqnitem 1.
First, we will give some definitions and briefly describe the most important known facts
concerning symmetric systems of polynomial equations,
which we will need later for presenting of the original part of the article.

A polynomial $P_s(x,y)$ in two variables is called {\it symmetric}
if it's value is unaltered  when the variables are interchanged:
$\,P_s(x,y)=P_s(y,x)$. In terms of transformations,
a symmetric polynomial is defined as a polynomial that
does not change its form under the transformation $x=\bar y$, $y=\bar x$.

The simplest symmetric polynomials
\begin{align}
\sigma_1=x+y,\quad \sigma_2=xy
\label{e01}
\end{align}
are called {\it elementary} polynomials. These polynomials are the
simplest algebraic invariants that do not change when the variables are permuted.

Let us denote
\begin{equation}
s_n=x^n+y^n,
\label{e02}
\end{equation}
where $n$ is a positive integer.
\goodbreak

The recurrent formula is valid
\begin{align*}
s_k=\sigma_1s_{k-1}-\sigma_2s_{k-2},
\end{align*}
which makes it possible to consistently express the sums of powers of the form \eqref{e02}
through elementary polynomials. For example,the first six sums of powers $s_n$ ($n=1, \ldots, 6$)
are presented in Table~1.

\begin{table}[h!]
\renewcommand*{\arraystretch}{1.4}

\noindent \caption{{\bf Table 1:} Sums of powers $s_n = x^n + y^n$ in terms of $\sigma_1$ and $\sigma_2$
($\sigma_1 = x+ y$, $\sigma_2 = xy$)}

\centering
   \begin{tabular}{ll}
         \hline
				\qquad $s_n$ \qquad & \qquad Expression  in terms of $\sigma_1$, $\sigma_2$\\
			\hline
        \qquad $s_1$ \qquad &\qquad $\sigma_1$ \\
				\qquad $s_2$ \qquad &\qquad $\sigma1^2 - 2\sigma_2$ \\
			  \qquad $s_3$ \qquad &\qquad $\sigma_1^3 - 3\sigma_1\sigma_2$ \\
			  \qquad $s_4$ \qquad &\qquad $\sigma_1^4 - 4\sigma_1^2\sigma_2 + 2\sigma_2^2$ \\
				\qquad $s_5$ \qquad &\qquad $\sigma_1^5 - 5\sigma_1^3\sigma_2 + 5\sigma_1\sigma_2^2$ \\
			  \qquad $s_6$ \qquad &\qquad $\sigma_1^6 - 6\sigma_1^4\sigma_2 + 9\sigma_1^2\sigma_2^2 - 2\sigma_2^3$ \\
				\qquad $s_7$ \qquad &\qquad $\sigma_1^7 - 7\sigma_1^5\sigma_2 + 14\sigma_1^3\sigma_2^2 - 7\sigma_1\sigma_2^3$ \\
				\qquad $s_8$ \qquad &\qquad $\sigma_1^8 - 8\sigma_1^6\sigma_2 + 20\sigma_1^4\sigma_2^2 - 16\sigma_1^2\sigma_2^3 + 2\sigma_2^4$ \\
        \hline
    \end{tabular}
\end{table}
\medskip

However, Waring's formula \cite{bol2002,Dickson1922} is more convenient:
\begin{align*}
\frac 1ks_k=\frac 1k\sigma_1^k-\frac{(k-2)!}{1!\,(k-2)!}\sigma_1^{k-2}\sigma_2+\frac{(k- 3)!}{2!\,(k-4)!}\sigma_1^{k-4}\sigma_2^2-
\frac{(k-4)!}{3!\,(k-6)!}\sigma_1^{k-6}\sigma_2^3+\cdots\,
\end{align*}
which allows one to immediately obtain in explicit form the required representation
of the sum of powers.

{\sl {\bf Fundamental theorem on symmetric polynomials.}
Any symmetric polynomial in two variables can be uniquely
expressed in terms of elementary polynomials~\eqref{e01}}.

The proof of this theorem can be found, for example, in \cite{bol2002,Blum-Smith2017}.

For solving systems of two polynomial equations
\begin{equation}
P_s(x,y)=0,\quad \ Q_s(x,y)=0,
\label{e02a}
\end{equation}
where $P_s$ and $Q_s$ are symmetric polynomials, it is useful
to  apply elementary symmetric polynomials \eqref{e01} as new variables,
which leads to a simplification of the original system.
System \eqref{e02a} does not change when the variables $x\rightleftharpoons y$ are permuted.
Therefore, such systems have the following property:
if a system has a solution $x=x_0$, $y=y_0$,
it also has a solution $x=y_0$, $y=x_0$.

If $\sigma_1$ and $\sigma_2$ are solutions to the transformed system,
then the required quantities $x$ and $y$ are determined from
simple system \eqref{e01},
whose solution reduces to a quadratic equation for~$x$
and a linear relation for~$y$:
\begin{equation}
x^2-\sigma_1x+\sigma_2=0,\quad \ \ y=\sigma_1-x.
\label{e03}
\end{equation}

Systems of two equations of the form \eqref{e02a} and
related polynomial systems with three or more unknowns,
the equations of which do not change with any permutation of unknowns,
will be called \textit{classical symmetric systems of polynomial equations}
(or briefly \textit{symmetric systems of polynomial equations}).

\begin{figure}[htb]
\centering
\includegraphics[width=0.90\textwidth]{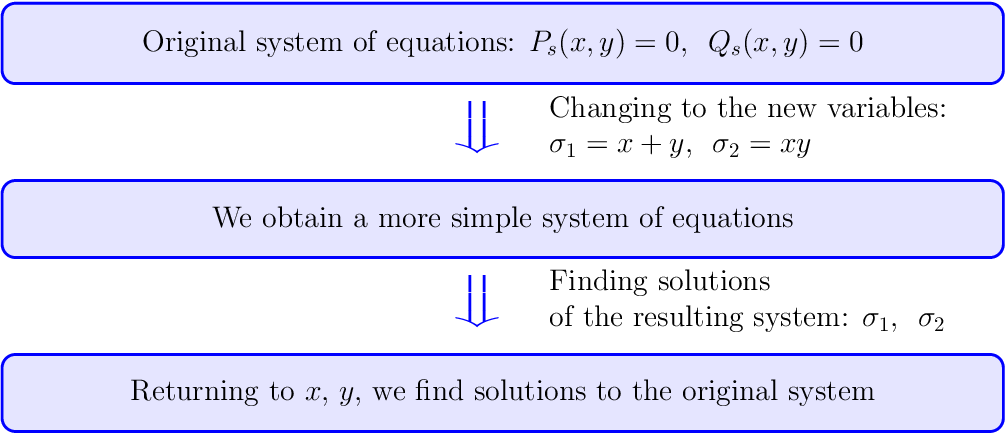}
\bigskip
\caption{{\bf Figure~1.} Scheme for solving classical symmetric systems of polynomial equations.}
\label{fig:1}
\end{figure}

For clarity, Fig.~\ref{fig:1} shows a schematic diagram for
solving symmetric systems of polynomial equations based on the use of
invariants \eqref{e01} as new variables.

\smallskip

\eqnitem 2.
Specific polynomial equations have ``hidden" (non-obvious, implicit) symmetries and,
by introducing a new additional variable, can be reduced to symmetric systems of polynomial equations.
Let us illustrate this with a specific example.
\medskip

\textit{Example 1.}
Consider the sixth-degree polynomial equation
\begin{equation}
(a-x^2)^3=(b-x^3)^2,
\label{e06}
\end{equation}
containing two free parameters $a$ and $b$ and in expanded form can be written as
\begin{equation}
2x^6-3ax^4-2bx^3+3a^2x^2+b^2-a^3=0.
\label{e07}
\end{equation}

Equation \eqref{e06} and the equivalent equation \eqref{e07}
do not have obvious symmetries.
Let us show that solutions to equation \eqref{e06} can be expressed
in terms of solutions to the symmetric polynomial system
\begin{equation}
\begin{aligned}
x^2+y^2&=a,\\
x^3+y^3&=b.
\end{aligned}
\label{e08}
\end{equation}
Indeed, the direct elimination of the additional variable~$y$
from system \eqref{e08} leads to equation~\eqref{e06}.

System \eqref{e08} is solved by changing to the new variables \eqref{e01}.
Considering the relations
\begin{align*}
x^2+y^2&=(x+y)^2-2xy=\sigma_1^2-2\sigma_2,\\
x^3+y^3&=(x+y)^3-3xy(x+y)=\sigma_1^3-3\sigma_1\sigma_2,
\end{align*}
we get a simpler system
\begin{align*}
\sigma_1^2-2\sigma_2&=a,\\
\sigma_1^3-3\sigma_1\sigma_2&=b.
\end{align*}
Eliminating $\sigma_2$, we arrive at the cubic equation
\begin{equation}
\begin{aligned}
\sigma_1^3-3a\sigma_1+2b=0.
\end{aligned}
\label{e09}
\end{equation}
Since the roots of any cubic equation can be expressed
in radicals~(see, for example, \text{\cite{bro2015,kor2000}}),
thus, the following statement has been constructively proven:
all the roots of the two-parameter sixth-degree equation \eqref{e06} can be expressed in radicals.

Considering the two degenerated cases, when $a=0$, $b\ne 0$ and $b=0$, $a\ne 0$,
the original polynomial equation \eqref{e06} is reduced, correspondingly, to the following
more simple polynomial equations:
\begin{equation}
\begin{aligned}
& 2x^6 - 2bx^3 + b^2 = 0,\quad {\rm or} \quad x^6+(b-x^3)^2 = 0\\
& 2x^6 - 3ax^4 + 3a^2x^2 - a^3 = 0,\quad {\rm or}  \quad x^6-(a-x^2)^3=0,
\end{aligned}
\label{eDeg}
\end{equation}
whose solutions are expressed in radicals.
\medskip

\textit{Remark 3.}
Instead of the sixth-degree polynomial equation \eqref{e06},
one can also consider more complicated higher-degree polynomial equations of the form
\begin{equation*}
(a-x^k)^n=(b-x^n)^k,
\end{equation*}
where $k$ and $n$ are positive integers,
which are derived from symmetric polynomial systems
\begin{align*}
x^k+y^k=a,\quad \ x^n+y^n=b.
\end{align*}

\textit{Remark 4.}
Symmetric systems of polynomial equations with respect to three or more unknowns,
preserving their form under all possible permutations of variables,
and methods for solving them are discussed in the book \cite{bol2002}.
This book also shows how symmetric polynomials can be used
to prove inequalities and to solve some irrational equations.

\section{Symmetric systems of polynomial equations of mixed type}

\eqnitem 1.
A polynomial $Q_a(x,y)$ in two variables is said to be {\it anti-symmetric}
if it changes sign when the variables are swapped, i.e. $\,Q_a(x,y)=-Q_a(y,x)$.

Two main properties of anti-symmetric polynomials are (see, for example, \cite{bol2002}):
\begin{equation}
\begin{aligned}
1.&\qquad Q_a(x,x)=0,\\
2.&\qquad Q_a(x,y)=(x-y)R_s(x,y),
\end{aligned}
\label{e10}
\end{equation}
where $R_s(x,y)$ is a symmetric polynomial.

\eqnitem 2. Consider a symmetric system of polynomial equations of mixed type
\begin{equation}
\begin{aligned}
P_s(x,y)&=0,\\
Q_a(x,y)&=0,
\end{aligned}
\label{e11}
\end{equation}
in which the first polynomial $P_s(x,y)$ is symmetric,
and the second polynomial $Q_a(x,y)$ is anti-symmetric.
When swapping the variables $x\rightleftharpoons y$,
the first equation of system \eqref{e11} does not change,
but the left-hand side of the second equation changes the sign
(and after multiplication by $-1$, this equation converts
into the original equation).
Therefore, such systems have the following property:
if system \eqref{e11} has a solution $x=x_0$, $y=y_0$,
it also has a solution $x=y_0$, $y=x_0$.

Due to the second property \eqref{e10}, system \eqref{e11}
can be represented in the form
\begin{equation}
\begin{aligned}
P_s(x,y)&=0,\\
(x-y)R_s(x,y)&=0.
\end{aligned}
\label{e12}
\end{equation}
The second equation in \eqref{e12} can be satisfied
by setting any of the two factors on its left-hand side equal to zero.
Therefore, system \eqref{e11} splits into two simpler subsystems
\begin{equation}
\begin{aligned}
P_s(x,x)&=0,\quad \ y=x;\\
P_s(x,y)&=0,\quad \ R_s(x,y)=0.
\end{aligned}
\label{e13}
\end{equation}
The first subsystem \eqref{e13} is reduced to a single polynomial equation.
The second subsystem \eqref{e13} is symmetric and
can be solved by the method described in Section~\ref{ss:1.2.4}.
\medskip

\textit{Example 2.}
Consider a two-parameter mixed polynomial system of the form \eqref{e11}:
\begin{equation}
\begin{aligned}
x^2+y^2&=a,\\
x^3+by&=y^3+bx.
\end{aligned}
\label{e14}
\end{equation}
After moving all terms to the left-hand side of the equations in~\eqref{e14},
we obtain a special case of system \eqref{e11} with
$$
P_s(x,y)=x^2+y^2-a,\quad \
Q_a(x,y)=x^3-y^3+b(y-x).
$$
The second equation \eqref{e14} can be represented as a product
\begin{equation*}
(x-y)(x^2+xy+y^2-b)=0.
\end{equation*}
Therefore, the original system \eqref{e14} splits into two simple subsystems
\begin{align*}
2x^2-a&=0,\quad \ y=x;\\
x^2+y^2-a&=0,\quad x^2+xy+y^2-b=0,
\end{align*}
the solution to which is elementary.
\medskip

\textit{Remark 5.}
In some cases, a system of two polynomial equations with two unknowns,
consisting of non-symmetric equations can be reduced to a symmetric system
of the form \eqref{e02a} by introducing new unknowns and applying a scaling transformation
of the form $x=\mu \bar x$, $y=\lambda \bar y$,
where $\mu$ and $\lambda$ are the required parameters.

\section{Nonclassical symmetric systems of polynomial equations}\label{sec4}

Let us now consider a system of polynomial equations of a special form
\begin{equation}
\begin{aligned}
P(x,y)+Q_s(x,y)&=0,\\
P(y,x)+Q_s(x,y)&=0,
\end{aligned}
\label{e15}
\end{equation}
where $P(x,y)$ is some given polynomial in two variables,
and $Q_s(x,y)$ is a symmetric polynomial.
When the variables are swapped, equations \eqref{e15} are interchangeable.
We will call such systems \textit{nonclassical symmetric systems of polynomial equations}.

If system \eqref{e15} has a solution
$x=x_0$, $y=y_0$, then it also has a solution $x=y_0$, $y=x_0$.

The solution of the system of equations \eqref{e15} can be reduced
to the solution of a more simple symmetric system and the solution
of a single independent equation.
To show this, first, by adding term by term to both equations,
we obtain the equation
\begin{equation}
P(x,y)+P(y,x)+2Q_s(x,y)=0,
\label{e16}
\end{equation}
which does not change when the variables are swapped.
Then subtracting the second equation from the first equation \eqref{e15},
we have
\begin{equation}
P(x,y)-P(y,x)=0.
\label{e17}
\end{equation}
It is easy to verify that the left-hand side of equation~\eqref{e17}
contains an anti-symmetric polynomial. Due to the second property~\eqref{e10},
the following representation is valid:
\begin{equation}
P(x,y)-P(y,x)=(x-y)R_s(x,y),
\label{e18}
\end{equation}
where $R_s(x,y)$ is a symmetric polynomial in two variables.
Therefore, equation \eqref{e17} can be written as
\begin{equation}
(x-y)R_s(x,y)=0.
\label{e19}
\end{equation}
Thus, the original system \eqref{e15} is equivalent to a system
consisting of two equations \eqref{e16} and \eqref{e19}.
Equation \eqref{e19} can be satisfied by setting any of the two factors
on its left-hand side equal to zero.
As a result, system \eqref{e16} and \eqref{e19}
splits into two more simple subsystems
\begin{equation}
\begin{aligned}
P(x,y)+P(y,x)+2Q_s(x,y)&=0,\quad \ y=x;\\
P(x,y)+P(y,x)+2Q_s(x,y)&=0,\quad \ R_s(x,y)=0.
\end{aligned}
\label{e20}
\end{equation}

The first subsystem \eqref{e20} is simplified
after substituting $y=x$ from the second equation and can be written as
$$
P(x,x)+Q_s(x,x)=0,\quad \ y=x.
$$
The second subsystem \eqref{e20} is symmetric and can be solved
by the method described in Section~\ref{ss:1.2.4}.

\begin{figure}[htb]
\centering
\includegraphics[width=0.90\textwidth]{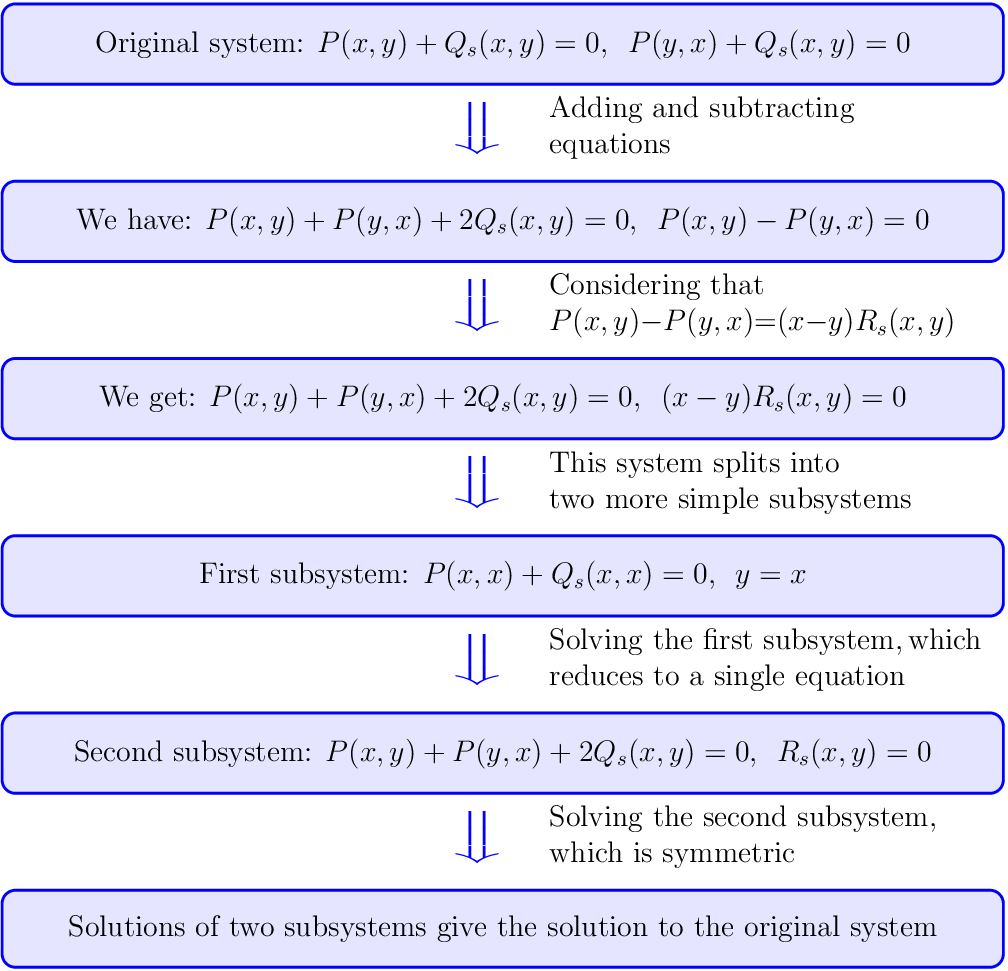}
\bigskip
\caption{{\bf Figure~2.} Scheme for solving nonclassical symmetric systems of polynomial equations.}
\label{fig:2}
\end{figure}

For clarity, Fig.~\ref{fig:2} shows a schematic diagram
for solving nonclassical symmetric systems of polynomial equations.

\smallskip
\bigskip

\medskip

\textit{Example 3.}
Consider a three-parameter polynomial system of the form \eqref{e15}:
\begin{equation}
\begin{aligned}
x&=ax^2+by^2+c,\\
y&=ay^2+bx^2+c,
\end{aligned}
\label{e21}
\end{equation}
where $a\not=\pm b$.
After moving all terms to the right-hand side of the equations
in~\eqref{e21}, we obtain a special case of system \eqref{e15} with
$$
P(x,y)=ax^2-x+by^2+c,\quad \ Q_s(x,y)\equiv 0.
$$

By adding and subtracting term by term the equations of system~\eqref{e21},
we arrive at the equivalent system of equations
\begin{align*}
(a+b)(x^2+y^2)-(x+y)+2c&=0,\\
(x-y)[(a-b)(x+y)-1]&=0,
\end{align*}
which splits into two simple subsystems
\begin{align*}
(a+b)x^2-x+c&=0,\quad \ y=x;\\
(a+b)(x^2+y^2)-(x+y)&=0,\quad \ (a-b)(x+y)-1=0.
\end{align*}
The solution of these subsystems is omitted due to its simplicity.

\section{Polynomial equations reducible to nonclassical\\ symmetric systems}

\eqnitem 1. Polynomial equations of the form
\begin{equation}
f(f(x))=x,
\label{e22}
\end{equation}
are reduced to the systems of type \eqref{e15} for $Q_s(x,y)\equiv 0$,
on the left-hand side of which there is the second iteration of
the given polynomial $f(x)$.
The roots of the more simple equation $f(x)=x$ are also the roots of
the original equation \eqref{e22} (this property is also valid
for more complicated transcendental equations, when $f(x)$
is any given continuous function).

The introduction of a new variable $y=f(x)$ allows us to obtain
from equation~\eqref{e22} the following equivalent polynomial system:
\begin{equation}
y=f(x),\quad \ x=f(y),
\label{e23}
\end{equation}
which is a special case of system \eqref{e15}
(after moving all terms of equations~\eqref{e23} to the left-hand side)
and is solved by the method described in Section~\ref{sec4}.
\medskip

\textit{Example 4.}
Consider the fourth-degree polynomial equation
\begin{equation}
(x^2+a)^2+a=x,
\label{e24}
\end{equation}
containing a free parameter $a$ and can be represented
as an equation~\eqref{e22} with $f(x)=x^2+a$.

By introducing a new variable $y=x^2+a$, equation \eqref{e24} is reduced to
a nonclassical symmetric system of the form \eqref{e15}:
\begin{equation}
y=x^2+a,\quad \
x=y^2+a.
\label{e25}
\end{equation}
Adding and subtracting the equations \eqref{e25},
as described in Section~\ref{sec4}, and after elementary transformations
and factorization of the second reduced equation, we have
\begin{equation}
\begin{aligned}
x^2+y^2-x-y+2a&=0,\\
(y-x)(x+y+1)&=0.
\end{aligned}
\label{e26}
\end{equation}
Further, equating to zero the factors on the left-hand side of the
second equation \eqref{e26}, in the end we obtain
the two independent quadratic equations
\begin{align*}
x^2-x+a&=0,\quad \ y=x;\\
x^2+x+a+1&=0,\quad \ y=-x-1,
\end{align*}
which determine the four roots $x$ of the original equation \eqref{e24}.
\medskip

\textit{Example 5.}
Let us consider a more complicated ninth-degree polynomial equation
of the form~\eqref{e22} with $f(x)=x^3+a$:
\begin{equation}
(x^3+a)^3+a=x.
\label{e27}
\end{equation}
This equation contains a free parameter $a$ and reduces to a nonclassical symmetric system
\begin{equation}
y=x^3+a,\quad \
x=y^3+a.
\label{e28}
\end{equation}
By adding and subtracting the equations \eqref{e28},
as described in Section~\ref{sec4}, and after elementary transformations,
we arrive at the equivalent system of equations
\begin{align*}
x^3+y^3-x-y+2a&=0,\\
(y-x)(x^2+xy+y^2+1)&=0,
\end{align*}
which splits into two more simple subsystems
\begin{equation}
\begin{aligned}
x^3-x+a&=0,\quad \ y=x;\\
x^3+y^3-x-y+2a&=0,\quad \ x^2+xy+y^2+1=0.
\end{aligned}
\label{e29}
\end{equation}
The solutions of the first subsystem \eqref{e29} are determined
by the roots of the cubic equation, which allows one to find three roots
of the original equation \eqref{e27}.
To solve the second system \eqref{e29}, which is symmetric and
determines the remaining six roots of the equation \eqref{e27},
we introduce the new variables by using formulas~\eqref{e01}.
As a result, we obtain the more simple system
\begin{equation}
\sigma_1^3-3\sigma_1\sigma_2-\sigma_1+2a=0,\quad \
\sigma_1^2-\sigma_2+1=0.
\label{e30}
\end{equation}
By eliminating $\sigma_2$ from the last system,
we arrive at the cubic equation for $\sigma_1$:
\begin{equation}
\sigma_1^3+2\sigma_1-a=0.
\label{e31}
\end{equation}
The roots of this equation and the relation $\sigma_2=\sigma_1^2+1$,
which follows from the second equation in~\eqref{e30},
give three pairs of real or complex numbers $\sigma_{1k}$, $\sigma_{2k}$ ($k=1,\,2,\,3$).
By substituting these numbers into the quadratic equation \eqref{e03},
we can find the six roots of the original equation~\eqref{e27}.

Since the roots of any cubic equation can be expressed in radicals,
therefore, the following statement is constructively proven:
all the roots of one-parameter ninth-degree equation~\eqref{e27} can be expressed in radicals.

It should be noted that the original polynomial equation \eqref{e27} can be represented
in the factorized form
\begin{equation}
(x^3 + a - x)(x^6 + 2ax^3 + x^4 + a^2 + ax + x^2 + 1)=0.
\label{ee27}
\end{equation}
Therefore, the solution of the nine-degree equation~\eqref{e27} consists of the three solutions
of the one-parameter cubic equation   $x^3+a-x=0$
and the six solutions of the one-parameter six-degree equation
$x^6 + 2ax^3 + x^4 + a^2 + ax + x^2 + 1=0$.
All nine solutions are expressed in radicals.
\medskip

\textit{Remark 6.}
If we are only interested in the real roots of the polynomial equation~\eqref{e27},
then the second system \eqref{e29} must be discarded, since in this case the
following inequality:
\begin{equation*}
x^2+xy+y^2+1=(x+\tfrac12y)^2+\tfrac34y^2+1>0
\end{equation*}
is valid.
\medskip

\textit{Remark 7.}
Similarly, it can be shown that a four-parameter polynomial equation of the ninth degree of the form \eqref{e22}
with $f(x)=a_3x^3+a_2x^2+a_1x+a_0$, which generalizes equation \eqref{e27}, is also solvable in radicals.
In particular, by setting $f(x)=(x+a)^3$, one can obtain a one-parameter equation solvable in radicals, $[(x+a)^3+a]^3-x=0$.

\eqnitem 2. Let us show that polynomial equations of the form
\begin{equation}
f(af(x)+x+ab)+f(x)+2b=0,
\label{e32}
\end{equation}
where $f(x)$ is a given polynomial, $a$ and $b$ are free parameters,
have ``hidden" symmetry and are reduced to the systems of type \eqref{e15}
for $Q_s(x,y)\equiv 0$.
Note that the roots of the more simple equation $f(x)=-b$ are also the roots of the
original equation \eqref{e32} (this property is also valid for more
complicated transcendental equations, when $f(x)$ is any given continuous function).

Let us denote
\begin{equation}
y=af(x)+x+ab.
\label{e33}
\end{equation}
As a result, equation \eqref{e32} takes the form
\begin{equation}
f(y)+f(x)+2b=0.
\label{e34}
\end{equation}
By eliminating the function $f(x)$ from \eqref{e33} and \eqref{e34}, we obtain the relation
\begin{equation}
x=af(y)+y+ab.
\label{e35}
\end{equation}
Equations \eqref{e33} and \eqref{e35} form a system of the form \eqref{e15}:
\begin{equation}
y=af(x)+x+ab,\quad \ x=af(y)+y+ab,
\label{e36}
\end{equation}
that can be solved by the method described in Section~\ref{sec4}.
\medskip

\textit{Example 6.}
Consider the one-parameter fourth-degree polynomial equation
\begin{equation}
(x^2+x+b)^2+x^2+2b=0,
\label{e37}
\end{equation}
which is a special case of equation \eqref{e32} (for $f(x)=x^2$, $a=1$)
and reduces to a nonclassical symmetric system
\begin{equation}
y=x^2+x+b,\quad \ x=y^2+y+b.
\label{e38}
\end{equation}
By adding and subtracting the equations \eqref{e38}, we obtain the equivalent system
of equations, which can be represented in the form
\begin{equation}
x^2+y^2+2b=0,\quad \
(y-x)(x+y+2)=0.
\label{e39}
\end{equation}
By equating to zero the factors on the left-hand side of the second equation, and
after elementary transformations, in the end, we arrive at
the two independent quadratic equations for the variable $x$
\begin{align*}
x^2+b&=0,\quad \ y=x;\\
x^2+2x+b+2&=0,\quad \ y=-x-2.
\end{align*}
These equations have roots
$$
x_1=-\sqrt{-b},\quad \ x_2=\sqrt{-b},\quad x_3=-1-\sqrt{-b-1},\quad x_4=-1+\sqrt{-b -1},
$$
which determine the required solutions to the original equation \eqref{e39}.
\medskip

\textit{Example 7.}
It can be shown that the one-parameter ninth-degree polynomial equation
\begin{equation}
(x^3+x+b)^3+x^3+2b=0,
\label{e32ab}
\end{equation}
which is a special case of equation \eqref{e32} (with $f(x)=x^3$, $a=1$),
is solvable in radicals.

It should be noted that the polynomial equation \eqref{e32ab} can be represented
in the factorized form
\begin{equation}
(x^3 + b)(x^6 + 2bx^3 + 3x^4 + b^2 + 3bx + 3x^2 + 2) = 0.
\label{e32ac}
\end{equation}
Therefore, the solution of the nine-degree equation~\eqref{e32ab} consists of the three solutions
of the one-parameter cubic equation   $x^3+b=0$
and the six solutions of the one-parameter six-degree equation
$x^6 + 2bx^3 + 3x^4 + b^2 + 3bx + 3x^2 + 2 = 0$.
All nine solutions are expressed in radicals.

\section{Other polynomial systems splitting into simpler\\ subsystems}

An important property of the system \eqref{e15} is that its equations
coincide for $y=x$. Let us consider a more general polynomial system of two equations
\begin{equation}
P(x,y)=0,\quad \ Q(x,y)=0,
\label{eee01}
\end{equation}
where polynomials in two variables $P(x,y)$ and $Q(x,y)$ for some constants
$\lambda$ and $\mu$ satisfy the condition
\begin{equation}
P(x,\lambda x)\equiv \mu Q(x,\lambda x).
\label{eee02}
\end{equation}
From \eqref{eee02}, it follows that the following representation
\begin{equation}
P(x,y)-\mu Q(x,y)=(y-\lambda x)R(x,y)
\label{eee03}
\end{equation}
takes place, where $R(x,y)$ is some polynomial.

Now, instead of the system \eqref{eee01}, we consider an equivalent system,
leaving the first equation and replacing the second equation with a
linear combination of equations $\,P(x,y)-\mu Q(x,y)=0$.
Taking into account the representation \eqref{eee03}, we obtain the system
\begin{equation}
P(x,y)=0,\quad \ (y-\lambda x)R(x,y)=0,
\label{eee04}
\end{equation}
which splits into two more simple subsystems
\begin{equation}
\begin{aligned}
P(x,y)&=0,\quad \ y=\lambda x;\\
P(x,y)&=0,\quad \ R(x,y)=0.
\end{aligned}
\label{eee05}
\end{equation}

\begin{figure}[h!tb]
\centering
\includegraphics[width=0.90\textwidth]{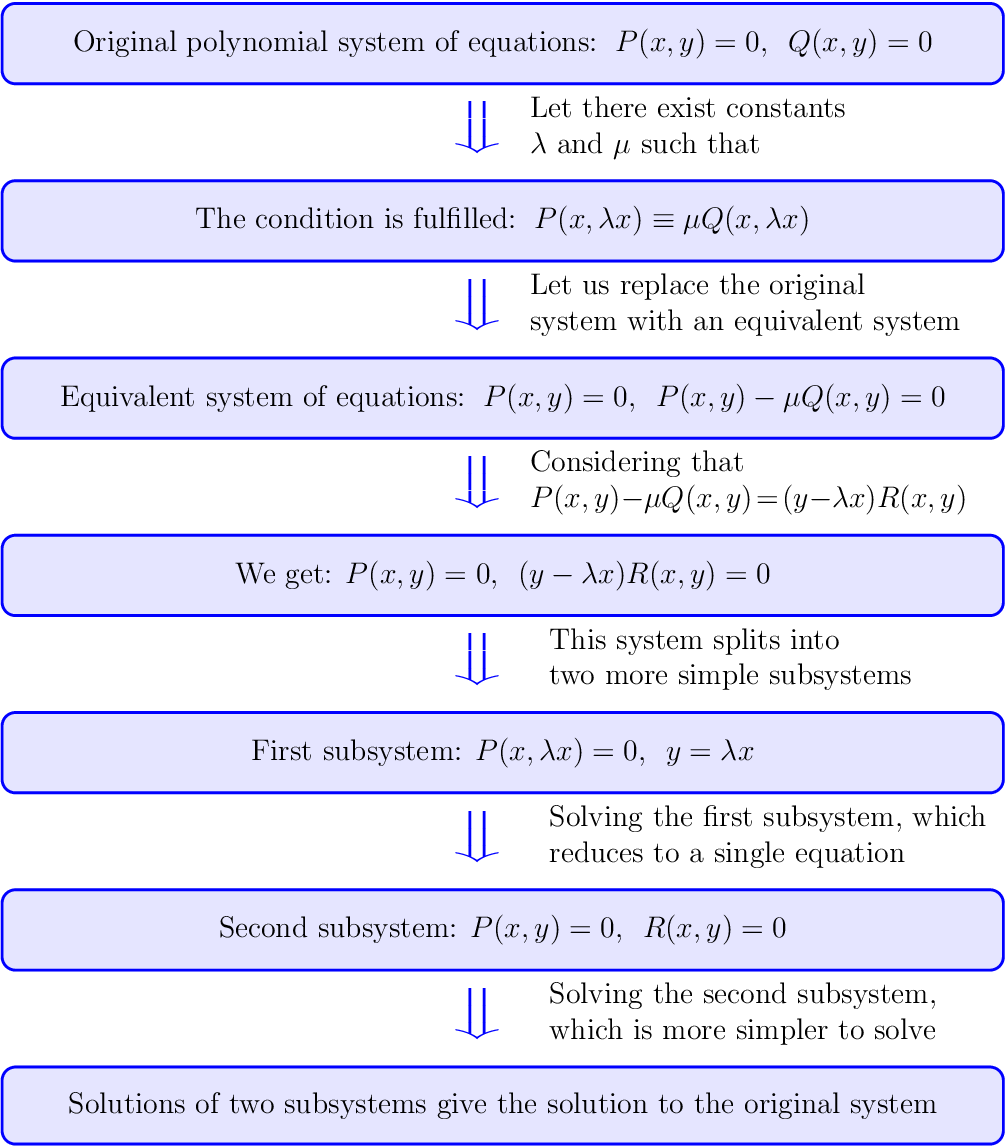}
\bigskip
\caption{{\bf Figure~3.} Scheme for solving systems of polynomial equations of the form \eqref{eee01} under condition \eqref{eee02}.}
\label{fig:3}
\end{figure}

For clarity, Fig.~\ref{fig:3} shows a schematic diagram
for solving systems of polynomial equations of the form \eqref{eee01} under condition \eqref{eee02}.
\medskip

\textit{Example 8.} Consider the polynomial system
\begin{equation}
ax^2+by^2+x+c=0,\quad \ (a+b)x^2-y+c=0,
\label{eee06}
\end{equation}
which is a special case of a system of the form \eqref{eee01},
which satisfies the condition \eqref{eee02} with $\lambda =-1$ and $\mu=1$.

Leaving the first equation and replacing the second equation with
the difference of equations, we obtain an equivalent system,
which can be represented as
\begin{equation*}
ax^2+by^2+x+c=0,\quad \ (x+y)[b(y-x)+1]=0.
\end{equation*}
This system splits into two simple subsystems
\begin{align*}
&ax^2+by^2+x+c=0,\quad \ y=-x;\\
&ax^2+by^2+x+c=0,\quad \ b(y-x)+1=0,
\end{align*}
each of which, after eliminating $y$, reduces to a quadratic equation for $x$.

\bigskip

\section{Test problems with parameters for analyzing\\ the capabilities of computer algebra systems}\label{sec7}

In this section, we formulate three test problems with parameters based on the previously considered
high-degree polynomial equations (containing free parameters) whose roots are expressed in radicals
(through solutions of cubic equations).
These problems will be used in subsequent sections to analyze and compare
the capabilities of the main computer algebra systems.

\medskip
\textbf{Test problem 1.}
As the first test problem, we will take a sixth-degree polynomial equation with two parameters:
\begin{equation}
(a-x^2)^3=(b-x^3)^2,
\label{ee06}
\end{equation}
all roots of which can be expressed in radicals (for more details, see Example~1).

\medskip
\textbf{Test problem 2.}
As the second test problem, we will take a more complicated
ninth-degree polynomial equation that contains one parameter:
\begin{equation}
(x^3+a)^3+a=x,
\label{eee27}
\end{equation}
all roots of which can also be expressed in radicals
(for more details, see Example~5).

\medskip
\textbf{Test problem 3.}
As the third test problem, we will take another ninth-degree polynomial equation
with one parameter:
\begin{equation}
(x^3+x+b)^3+x^3+2b=0,
\label{eeee27}
\end{equation}
all roots of which can be expressed in radicals (see Example~7).
\medskip


\section{Testing Maple}

In this section, we consider the test problems with parameters formulated in Section~\ref{sec7}
for analyzing the capabilities of one of the leading computer algebra systems Maple
(Maple 2024 version), for the analytical and numerical solution of polynomial equations.
\medskip

\textbf{Test problem 1.}
\noindent Now let us try to solve the six-degree polynomial equation \eqref{ee06} with Maple:
\begin{verbatim}
   Eq1:= (a-x^2)^3=(b-x^3)^2;
   S11:= solve(Eq1, x);
   S12:= allvalues(S11,explicit);
   S13:= convert([allvalues(S11,explicit)],radical);
\end{verbatim}
If $a$ and $b$ are arbitrary parameters, the results cannot be obtained explicitly.
In this case, Maple gives solutions only in terms of the Maple predefined function
{\tt RootOf} as follows:\footnote{The Maple predefined function {\tt RootOf}  is a way of
representing the roots of an equation when explicit formulas are too cumbersome,
or Maple cannot find them or do not exist at all.}:
{\small
\begin{equation*}
   S11:= RootOf(2\,\_Z^6-3\,\_Z^4a-2\,\_Z^3b+3\,\_Z^2a^2-a^3+b^2),
\end{equation*}}
\!\!where the global variable $\_Z$ is the variable we want to solve for this equation,
which is $x$ in this example.
If there are several roots (i.e., six in this example), Maple gives the result:
{\small
\begin{align*}
& S12:= RootOf(2\,\_Z^6-3\,\_Z^4a-2\,\_Z^3b+3\,\_Z^2a^2-a^3+b^2, index=1),\\
& RootOf(2\,\_Z^6-3\,\_Z^4a-2\,\_Z^3b+3\,\_Z^2a^2-a^3+b^2, index=2), \\
& RootOf(2\,\_Z^6-3\,\_Z^4a-2\,\_Z^3b+3\,\_Z^2a^2-a^3+b^2, index=3), \\
& RootOf(2\,\_Z^6-3\,\_Z^4a-2\,\_Z^3b+3\,\_Z^2a^2-a^3+b^2, index=4), \\
& RootOf(2\,\_Z^6-3\,\_Z^4a-2\,\_Z^3b+3\,\_Z^2a^2-a^3+b^2, index=5), \\
& RootOf(2\,\_Z^6-3\,\_Z^4a-2\,\_Z^3b+3\,\_Z^2a^2-a^3+b^2, index=6),
\end{align*}}
\!\!where the  {\tt index = ...}  indicates the ordinal number of the root.
To get the roots explicitly, we can use the Maple predefined function {\tt allvalues}.
However, since equation \eqref{ee06} also contains two arbitrary parameters  $a$ and $b$,
to obtain all the roots in explicit form,
we need to specify the specific numerical values of these parameters.

\medskip
For the two special cases \eqref{ee06}, when $a=0$, $b\ne 0$ and $b=0$, $a\ne 0$,
\begin{verbatim}
  Eq12:= (a-x^2)^3-(b-x^3)^2=0;
   S14:= expand(subs(a=0,Eq12));
   S15:= simplify([solve(S14, x)]);
   S16:= expand(subs(b=0,Eq12));
   S17:= solve(S16, x);
\end{verbatim}
Maple gives all the six solutions expressed in radicals:
{\small
\begin{align*}
&\frac12 \left((4+4I)b\right)^{1/3}, \frac14 \left((4+4I)b\right)^{1/3}(I\sqrt{3}-1),
 -\frac14 \left((4+4I)b\right)^{1/3}(I\sqrt{3}+1), \\
&\frac12 \left((4-4I)b\right)^{1/3}, \frac14 \left((4-4I)b\right)^{1/3}(I\sqrt{3}-1),
  -\frac14 \left((4-4I)b\right)^{1/3}(I\sqrt{3}+1),
\end{align*}
}
and
{\small
\begin{align*}
& \frac12 \sqrt{2}\sqrt{a}, -\frac12 \sqrt{2}\sqrt{a}, \frac12\sqrt{2Ia\sqrt{3}+2a},
-\frac12 \sqrt{2Ia\sqrt{3}+2a}, \\
& \frac12\sqrt{-2Ia\sqrt{3}+2a}, -\frac12\sqrt{-2Ia\sqrt{3}+2a}.
\end{align*}
}

\textit{Remark 8.}
In both of the special cases discussed above, equation \eqref{ee06} is significantly simplified, namely,
for $a=0$, $b\ne 0$, the replacement $y=x^3$ leads it to a quadratic equation,
and for $b=0$, $a\ne 0$, the substitution $z=x^2$ leads to a cubic equation.
\medskip

Setting in the Maple script, various integer values of the parameters, for example, $a=5$ and $b=2$,
\begin{verbatim}
   allvalues(subs({a=5,b=2},S13));
\end{verbatim}
Maple gives all the six roots expressed in radicals:
\begin{equation}
\begin{aligned}
&  1-(1/2)\sqrt{3}+(1/2)\sqrt{3+4\sqrt{3}},\quad 1+(1/2)\sqrt{3}+(1/2)I\sqrt{4\sqrt{3}-3}, \\
& -2+(1/2 I)\sqrt{6},\quad 1-(1/2)\sqrt{3}-(1/2)\sqrt{3+4\sqrt{3}}, \\
& -2-(1/2 I)\sqrt{6},\quad 1+(1/2)\sqrt{3}-(1/2)I\sqrt{4\sqrt{3}-3},
\end{aligned}
\label{Roots3aa}
\end{equation}
where two roots are real. This means that in a certain range of parameter changes
(in the vicinity of these values, i.e. $a=5$, $b=2$) there will also be two real roots.

The other set of the parameters, for example, the integer values $a=7$ and $b=2$,
does not bring the results expressed in radicals.

If we specify the same parameters as the real numbers, $a=7.0$ and $b=2.0$,
\begin{verbatim}
   allvalues(subs({a=7.0,b=2.0}, S13));
\end{verbatim}
we get the corresponding numerical values (real and complex):
\begin{align*}
& 1.963798039, 2.242095980 + 1.235716141I, -2.337499474 + 1.401393518I,\\
& -1.772991050, -2.337499474 - 1.401393518I, 2.242095980 - 1.235716141I.
\end{align*}

\textit{Remark 9.}
Maple includes various specialized
predefined functions, packages, and sub-packages for classification, analysis, and finding
real solutions of polynomial equations with parameters,
see for example, package {\tt SolveTools}, predefined function {\tt SemiAlgebraic} with
option {\tt parameters}; package {\tt Root} {\tt Finding} with sub-package {\tt Parametric};
package {\tt RegularChains} with sub-packages {\tt ParametricSystemTools}
and {\tt SemiAlgebraicSetTools}.
However, these additional Maple options do not allow finding real roots in explicit form
of Test problem 1 with parameters. For example, although there are two real roots
in the segment $a=5$, $1\le b\le 3$, Maple does not give explicit expressions of these
real roots in this specified parameter domain,
nevertheless for $a=5$, $b=2$ these roots are expressed by formulas \eqref{Roots3aa}.

\medskip
\textbf{Test problem 2.}
Now let us try to solve the ninth-degree  polynomial equation \eqref{eee27} with Maple:
\begin{verbatim}
   Eq2:= (x^3+a)^3+a=x;
   S21:= [solve(Eq2, x)];
   S22:= allvalues(S21,explicit);
\end{verbatim}
As a result, we have the three roots expressed in radicals in explicit form
{\small
\begin{align}
\begin{split}
&\frac{1}{6}\left({}-108a {+} 12\sqrt{81a^2 {-} 12}\right)^{1/3} + 2\left({-}108a {+} 12\sqrt{81a^2 {-} 12}\right)^{-1/3},\\
&{-}\frac{1}{12}\left({-}108a {+} 12\sqrt{81a^2 {-} 12}\right)^{1/3} - \left({-}108a {+} 12\sqrt{81a^2 {-} 12}\right)^{-1/3} \\
&   \quad {+} \frac{I\sqrt{3}}{2}\left( \frac{1}{6}\left({-}108a {+} 12\sqrt{81a^2 {-} 12}\right)^{1/3}
          {-} 2\left({-}108a {+} 12\sqrt{81a^2 {-} 12}\right)^{-1/3}\right),\\
&{-}\frac{1}{12}\left({-}108a {+} 12\sqrt{81a^2 {-} 12}\right)^{1/3} {-} \left({-}108a {+} 12\sqrt{81a^2 {-} 12}\right)^{-1/3} \\
&   \quad {-} \frac{I\sqrt{3}}{2}\left( \frac{1}{6}\left(-108a + 12\sqrt{81a^2 {-} 12}\right)^{1/3}
          {-} 2\left({-}108a {+} 12\sqrt{81a^2 {-} 12}\right)^{-1/3}\right),	
\end{split}
\label{Roots3}
\end{align}
}
\!\!that is, these expressions are the roots of the cubic equation $x^3+a-x=0$ (see \eqref{ee27}),
and the other six roots of the six-degree equation $x^6 + 2ax^3 + x^4 + a^2 + ax + x^2 + 1=0$
are expressed in terms of the predefined function {\tt RootOf}.

\medskip
\noindent If we specify the integer value of the parameter~$a$, for example, $a=3$,
\begin{verbatim}
   S23:=allvalues(subs(a=3,[S22]));
   S23[1]; S23[2]; S23[3]; S23[4]; S23[5]; S23[6];
\end{verbatim}
we obtain all the nine roots of this equation expressed in radicals.
The expression {\tt S23} consists of six parts: in every part, the first three roots
are the same (that correspond to the explicit part of the result with a free parameter) and the fourth root is
different (that corresponds to the {\tt RootOf} part of the result with a free parameter).

\medskip
\noindent If we specify the same parameter as the real number, $a=3.0$,
\begin{verbatim}
   S24:=allvalues(subs(a=3.0,[S22]));
   S24[1]; S24[2]; S24[3]; S24[4]; S24[5]; S24[6];
\end{verbatim}
we get the corresponding numerical values (real and complex):
As a result, we get the expression that contains six parts:
{\small
\begin{align*}
& [0.835849940828641 {-} 1.04686931885012I, {-}1.67169988165728 {-} (2.\cdot 10^{{-}15})I, \\
& 0.835849940828640 {+} 1.04686931885012I, 0.844461199078011 {+} 1.39735098486455I] \\
& [0.835849940828641 {-} 1.04686931885012I, {-}1.67169988165728 {-} (2.\cdot 10^{{-}15})I, \\
& 0.835849940828640 {+} 1.04686931885012I, 0.500000000000000 {+} 1.32287565553230I] \\
& [0.835849940828641 {-} 1.04686931885012I, {-}1.67169988165728 {-} (2.\cdot 10^{{-}15})I, \\
& 0.835849940828640 {+} 1.04686931885012I, {-}1.34446119907801 {+} 0.260961410313150I] \\
& [0.835849940828641 {-} 1.04686931885012I, {-}1.67169988165728 {-} (2.\cdot 10^{{-}15})I, \\
& 0.835849940828640 {+} 1.04686931885012I, {-}1.34446119907801 {-} 0.260961410313150I] \\
& [0.835849940828641 {-} 1.04686931885012I, {-}1.67169988165728 {-} (2.\cdot 10^{{-}15})I, \\
& 0.835849940828640 {+} 1.04686931885012I, 0.500000000000000 {-} 1.32287565553230I] \\
& [0.835849940828641 {-} 1.04686931885012I, {-}1.67169988165728 {-} (2.\cdot 10^{{-}15})I, \\
& 0.835849940828640 {+} 1.04686931885012I, 0.844461199078011 {-} 1.39735098486455I].
\end{align*}
}
\!\!In every part, the first three roots are the same
(that correspond to the explicit part of the result with a free parameter)
and the fourth root is different (that corresponds to the {\tt RootOf} part
of the result with a free parameter).
\medskip

\textbf{Test problem 3.}
Now let us try to solve the ninth-degree polynomial equation \eqref{eeee27} with Maple:
\begin{verbatim}
   Eq3:= (x^3+x+b)^3+x^3+2*b=0;
   S31:= solve(Eq3, x);
   S32:= allvalues(S31,explicit);
\end{verbatim}
As a result, we have the three roots expressed in radicals in explicit form
\begin{align}
& (-b)^{1/3}, -\frac12(\sqrt{3}I + 1)(-b)^{1/3}, \frac12(-b)^{1/3}(\sqrt{3}I - 1),
\label{RootsTP3}
\end{align}
that is, these expressions are the roots of the cubic equation $x^3+b=0$ (see \eqref{e32ac}),
and the other six roots of the six-degree equation $x^6 + 2bx^3 + 3x^4 + b^2 + 3bx + 3x^2 + 2 = 0$
are expressed in terms of the predefined function {\tt RootOf}.

\medskip
\noindent If we specify the integer value of the parameter~$a$, for example, $a=4$,
we obtain all the nine roots of this equation expressed in radicals.
If we specify the same parameter as the real number, $a=4.0$,
we get the corresponding numerical values (real and complex) of the nine roots.

\section{Testing Mathematica}

In this section, we consider the test problems with parameters formulated in Section~\ref{sec7}
for analyzing the capabilities of one of the leading computer algebra systems Mathematica
(Mathematica 13 version) for the analytical and numerical solution of polynomial equations.
\medskip

\textbf{Test problem 1.}
If we solve the sixth-degree polynomial equation \eqref{ee06} with two parameters by using Mathematica,
\begin{verbatim}
   eq1 = (a-x^2)^3==(b-x^3)^2
   s11 = Solve[eq1,x]
   s12 = Roots[eq1,x]
   s13 = ToRadicals[s11]
\end{verbatim}
we also cannot obtain explicit roots.

Mathematica gives the six roots expressed
in terms of the predefined function {\tt Root} as follows:
{\small
\begin{verbatim}
{{x->Root[-a^3+b^2+3a^2#1^2-2b#1^3-3a#1^4+2#1^6&,1]},
 {x->Root[-a^3+b^2+3a^2#1^2-2b#1^3-3a#1^4+2#1^6&,2]},
 {x->Root[-a^3+b^2+3a^2#1^2-2b#1^3-3a#1^4+2#1^6&,3]},
 {x->Root[-a^3+b^2+3a^2#1^2-2b#1^3-3a#1^4+2#1^6&,4]},
 {x->Root[-a^3+b^2+3a^2#1^2-2b#1^3-3a#1^4+2#1^6&,5]},
 {x->Root[-a^3+b^2+3a^2#1^2-2b#1^3-3a#1^4+2#1^6&,6]}}.
\end{verbatim}
}

\medskip
\noindent For the two special cases \eqref{eDeg}, when $a=0$, $b\ne 0$ and $b=0$, $a\ne 0$,
\begin{verbatim}
   eq12 = -(a-x^2)^3+(b-x^3)^2==0
   s14  = eq12/.{a->0}
   s15  = Solve[s14,x]//FullSimplify
   s16  = eq12/.{b- 0}
   s17  = Solve[s16,x]
\end{verbatim}
Mathematica gives all the six solutions expressed in radicals:
{\small
\begin{verbatim}
{x->(1/2-I/2)^(1/3) b^(1/3)},
{x->((-1)^(7/12) b^(1/3))/2^(1/6)},
{x->(1/2+I/2)^(1/3) b^(1/3)},
{x->-(((1-I) b^(1/3))/2^(2/3))},
{x->-(((1+I)b^(1/3))/2^(2/3))},
{x->-(((-1)^(5/12) b^(1/3))/2^(1/6))},
\end{verbatim}
}
and
{\small
\begin{verbatim}
{x->-(Sqrt[a]/Sqrt[2])}, {x->Sqrt[a]/Sqrt[2]},
{x->-Sqrt[a/2-1/2ISqrt[3]a]}, {x->Sqrt[a/2-1/2ISqrt[3]a]},
{x->-Sqrt[a/2+1/2ISqrt[3]a]}, {x->Sqrt[a/2+1/2ISqrt[3]a]}.
\end{verbatim}
}

\medskip
\noindent
Setting in the Mathematica script, various integer values of the parameters, for example, $a=5$ and $b=2$,
\begin{verbatim}
   Solve[eq1,x]/.{a->5,b->2}
\end{verbatim}
Mathematica gives explicitly only the two roots expressed in radicals
\begin{verbatim}
   {x->1/2(-4-I Sqrt[6])}, {x->1/2(-4+I Sqrt[6])},
\end{verbatim}
and the other four roots are represented, as above, in terms of the predefined function {\tt Root}.
The other set of the parameters, for example, $a=7$ and $b=2$, as in Maple,
does not bring the results expressed in radicals.

\medskip
\noindent
If we specify the same parameters as the real numbers, $a=7.0$ and $b=2.0$
\begin{verbatim}
   Solve[eq1,x]/.{a->7.0,b->2.0}
\end{verbatim}
we get the corresponding numerical values (real and complex), the same as in Maple.
\medskip

\textbf{Test problem 2.}
If we solve the ninth-degree polynomial equation \eqref{eee27} with Mathematica,
\begin{verbatim}
   eq2 = (x^3+a)^3+a==x
   s21 = Solve[eq2, x]
\end{verbatim}
we have, as in Maple, the three roots expressed in radicals in the explicit form \eqref{Roots3},
and the six roots expressed in terms of the predefined function {\tt Root}.

\medskip
\noindent
If we specify the integer value of the parameter~$a$, for example, as in Maple, $a=3$
\begin{verbatim}
   Solve[eq2,x]/.{a->3}
\end{verbatim}
we obtain only the five roots of this equation expressed in radicals
in explicit form and the four roots are expressed in terms of the predefined function {\tt Root}.

\medskip
\noindent
If we specify the same parameter as the real number, $a=3.0$,
\begin{verbatim}
   Solve[eq2, x] /. {a -> 3.0}
\end{verbatim}
we get all the nine numerical values (real and complex).
\medskip

\textbf{Test problem 3.}
If we solve the ninth-degree  polynomial equation \eqref{eeee27} with Mathematica,
\begin{verbatim}
   eq3 = (x^3+x+b)^3+x^3+2*b==0
   s31 = Solve[eq3,x]
\end{verbatim}
we have, as in Maple, the three roots expressed in radicals in explicit form
{\small
\begin{verbatim}
{x->-b^(1/3)}, {x->(-1)^(1/3)b^(1/3)},
{x->-(-1)^(2/3)b^(1/3)}
\end{verbatim}
}
and the six roots expressed in terms of the predefined function {\tt Root}.

\medskip
\noindent
If we specify the integer value of the parameter~$a$, for example, as in Maple, $a=4$
\begin{verbatim}
   Solve[eq3,x]/.{a->4}
\end{verbatim}
we obtain only the five roots of this equation expressed in radicals
in explicit form and the four roots are expressed in terms of the predefined function {\tt Root}.

\begin{table}[htb]
\footnotesize

\noindent {\bf Table 2:} The results of testing Maple and Mathematica on three text problems with parameters
    \vspace{4mm}
\centering
\thickmuskip=3mu
\def\vruleskip{\noalign{\nointerlineskip}
  \omit& height 3.7pt &&&&&&&&&&\cr \noalign{\nointerlineskip}}
\def\rlineskip{\vruleskip\hline\vruleskip}
\def\\{\crcr}
\ialign to \textwidth{\strut#& \vrule#\tabskip=.1em plus1.8em&
  \hfil\ccbox{#}\hfil& \vrule#&
  \hfil\ccbox{#}\hfil& \vrule#&
  \hfil\lcbox{#}\hfil& \vrule#&
  \hfil\ccbox{#}\hfil& \vrule#&
  \hfil\ccbox{#}\hfil& \vrule#\tabskip=0pt \relax\cr
\hline \vruleskip
&& \strut No.&&
   \strut Problem statement&&
   \strut Problem parameters&&
   \strut Maple: root found&&
   \omit\strut\hfil Mathematica: root found\hfil&\cr
\vruleskip
\hline \omit& height.5mm&\omit&&\omit&&\omit&&\omit&&\omit&\cr
\hline
\vruleskip
&& 1&& $(a{-}x^2)^3{=}(b{-}x^3)^2$&&
   $a$, $b$ are any \\ $a=0$, $b$ is any \\ $a$ is any, $b=0$ \\ $a=2$, $b$ is any \\ $a$ is any, $b=2$ \\ $a=5$, $b=2$ \\
	 $a=7$, $b=2$ \\ $a=7.0$, $b=2.0$ &&
   \!\!0\\ 6 \\ 6 \\ \!\!0\\ \!\!0\\ 6 \\ 0 \\ 6 &&
   \!\!0\\ 6 \\ 6 \\ \!\!0\\ \!\!0\\ 2 \\ 0 \\ 6 &\cr
\rlineskip
&& 2&& $(x^3{+}a)^3{+}a{=}x$&&
   $a$ is any \\ $a=3$ \\ $a=3.0$ &&
   3 \\ 9 \\ 9 &&
   3 \\ 5 \\ 9 &\cr
\rlineskip
&& 3&& $(x^3{+}x{+}b)^3{+}x^3{+}2b{=}0$&&
   $b$ is any \\ $b=4$ \\ $b=4.0$ &&
   3 \\ 9 \\ 9 &&
   3 \\ 5 \\ 9 &\cr
   \vruleskip \hline}
\end{table}

\medskip
\noindent
If we specify the same parameter as the real number, $a=4.0$,
\begin{verbatim}
   Solve[eq3,x]/.{a->4.0}
\end{verbatim}
we get all the nine numerical values (real and complex).
\medskip

\medskip

For clarity, Table 2 shows the final results of using the Maple and Mathematica systems
for three text problems containing one or two free parameters.

\section{Conclusions}

This article discusses nonclassical symmetries and reductions
of polynomial equations and systems of polynomial equations.
It is shown that by introducing a new additional variable,
some polynomial equations having ``hidden" symmetries can be transformed
into classical symmetric systems of polynomial equations.

It is shown that nonclassical symmetric systems of two polynomial equations,
that change places when the unknowns are swapped,
can be reduced to simpler symmetric systems and a single independent equation.
It has been established that some equations of a special form,
containing the second iteration of a given polynomial,
can be transformed to nonclassical symmetric systems of polynomial equations.

Examples of non-trivial polynomial equations of the sixth and ninth degrees are given,
containing one or more free parameters that can be resolved in radicals.
These equations are investigated computationally by using computer algebra systems
Maple and Mathematica and can be proposed as test problems with parameters for
analytical and numerical methods of solving polynomial equations.

It is shown that neither Maple nor Mathematica can explicitly express in radicals
all or even some roots of a sixth-degree polynomial equation with two free parameters.
However, in the two degenerated cases (when one of the parameters is equal to zero),
both systems represent analytical solutions in radicals.
For some integer values of the parameters, Maple can obtain the six roots of the equation,
but Mathematica can find only four roots. For numerical values of the parameters,
both computer algebra systems can find all numerical roots of the equation.
For the one-parameter ninth-degree polynomial equation, Maple and Mathematica can find
only the three solutions expressed in radicals.
For some integer values of the parameter, Maple can obtain the nine roots of the equation,
but Mathematica can find only five roots. For numerical values of the parameter,
both computer algebra systems can find all numerical roots of the equation.
The test problems with parameters proposed in this work can be used for further improvements
of existing computer algebra systems.

\end{document}